\documentclass{amsart}
\usepackage{amssymb,graphicx,doublespace}
\newtheorem{theorem}{Theorem}[section]

\newtheorem{conjecture}[theorem]{Conjecture}

\theoremstyle{remark}

\newtheorem{example}[theorem]{Example}
\newcommand{\GVE}{\mathbf{G}=(V,E)}
\newcommand{\bfG}{\mathbf{G}}
\newcommand{\bfH}{\mathbf{H}}

\newcommand{\lgn}{\lg n}

\newcommand{\gn}{\mathcal{G}_{n,1/2}}
\newcommand{\gnp}{\mathcal{G}_{n,p}}
\newcommand{\cgQ}{\mathcal{P}}

\begin{document}
\title[GRAPH PEBBLING]%
 {A Note on Graph Pebbling}
\author[A. CZYGRINOW]{Andrzej Czygrinow}
\address{Department of Mathematics\\
         Arizona State University\\
         Tempe, Arizona~85287\\
         U.S.A.}
\email{andrzej@math.la.asu.edu}
\author[G. HURLBERT]{Glenn~Hurlbert}
\address{Department of Mathematics\\
         Arizona State University\\
         Tempe, Arizona~85287\\
         U.S.A.}
\email{hurlbert@math.la.asu.edu}
\author[H. A. KIERSTEAD]{H.~A.~Kierstead}
\address{Department of Mathematics\\
         Arizona State University\\
         Tempe, Arizona~85287\\
         U.S.A.}
\email{kierstead@ASU.edu}
\author[W. T. TROTTER]{William.~T.~Trotter}
\address{Department of Mathematics\\
         Arizona State University\\
         Tempe, Arizona~85287\\
         U.S.A.}
\email{trotter@ASU.edu}
\thanks{The research of the fourth author is supported in part by
 the Office of Naval Research.}
% =============================================================
%             Abstract and Maketitle
% =============================================================
\begin{abstract}
We say that a graph $\bfG$ is \textit{Class~$0$} if its pebbling
number is exactly equal to its number of vertices.
For a positive integer $d$, let $k(d)$ denote the least positive
integer so that every graph $\bfG$ with diameter
at most $d$ and connectivity at least $k(d)$ is Class 0.
The existence of the function $k$ was conjectured
by Clarke, Hochberg and Hurlbert, who showed that if the function
$k$ exists, then it must satisfy $k(d)=\Omega(2^d/d)$.  In this 
note, we show that $k$ exists and satisfies $k(d)=O(2^{2d})$.  
We also apply this result to improve the upper bound on the
random graph threshold of the Class 0 property.
\end{abstract}
\maketitle
\newpage
% =============================================================
\section{Introduction}\label{s:intro}
% =============================================================
Let $\mathbb{N}_0$ denote the non-negative integers.
When $\GVE$ is a finite graph, a function $\phi:V\rightarrow\mathbb{N}_0$
is called a \textit{pebbling}.  The quantity $\sum_{x\in V}\phi(x)$ is 
called the \textit{size} of $\phi$; the size of $\phi$ is just the total 
number of pebbles assigned to vertices.  In what follows, we 
consider a simple rule by which one pebbling is transformed into 
another: \quad  Choose a vertex $x$ to 
which at least two pebbles have been assigned.  Remove two 
pebbles from $x$ and add one pebble to an adjacent vertex $y$.
A pebbling obtained from $\phi$ by a sequence of such transformations 
is called a \textit{descendant} of $\phi$.

Given a vertex $x\in V$ and a pebbling $\phi$, we
say that $\phi$ \textit{pebbles} $x$ provided $\phi(x)>0$. Similarly, we
say that $\phi$ has the \textit{potential to pebble} $x$ provided
that $\phi$ or one of its descendants pebbles $x$.
The \textit{pebbling number} of a graph $\GVE$,
denoted $f(\bfG)$, is then the least $p$ so that for any pebbling 
$\phi$ of size $p$ and any vertex $x\in V$, $\phi$ has the 
potential to pebble $x$.

Clearly, the pebbling number of a graph is at least as large
as the number of vertices.  Following~\cite{CHH},
we say a graph $\GVE$ is \textit{Class~$0$} if
$f(\bfG)=|V|$.  Graphs which do not belong to
Class~$0$ are said to be \textit{Class~$1$}.

In~\cite{CHH}, Clarke, Hochberg and Hurlbert gave
the following conjecture.
% -------------------------------------------------------------
\begin{conjecture}\label{c:CHH}
For each $d\ge1$, there exists a least positive 
integer $k(d)$ so that all graphs of diameter~$d$ and 
connectivity at least $k(d)$ belong to Class~$0$. 
\end{conjecture}
% -------------------------------------------------------------
This conjecture holds trivially when $d=1$, since a graph of 
diameter~1 is a complete graph.  Note that a path on $3$ 
points shows that $k(2)\ge 2$.  However, as noted
in~\cite{CHH}, the following example shows that $k(2)\ge3$.
% -------------------------------------------------------------
\begin{example}\label{e:exam}
Label the vertices of a $6$-cycle as
$x_1,\dots,x_6$ so that $\{x_i,x_{i+1}\}$ is an edge for
$i=1,2,\dots,5$.  Of course, $\{x_1,x_6\}$ is also an edge.
Then let $\bfG_1$ be the graph formed by adding the edges 
$\{x_1,x_3\}$ and $\{x_3,x_5\}$.  Also, let $\bfG_2$ be
the graph formed from $\bfG_1$ by adding the edge $\{x_1,x_5\}$.  
Then $\bfG_1$ and $\bfG_2$ are $2$-connected and have diameter~$2$.  
However, the pebbling number of both graphs is~$7$.
\end{example}
% -------------------------------------------------------------
Clarke, Hochberg and Hurlbert showed that $k(2)=3$, and they 
characterized all $2$-connected graphs with diameter~$2$ which belong 
to Class~$1$.  The two graphs $\bfG_1$ and $\bfG_2$ constructed in
the preceding example are the only such graphs on $6$ or fewer 
vertices.  

Using a ``blow-up'' of a path, Clarke, Hochberg and Hurlbert also
showed that if the function $k(d)$ exists, then it must 
satisfy $k(d)=\Omega(2^d/d)$.
% =============================================================
\section{The Principal Result}\label{s:proof}
% =============================================================
In this section, we settle Conjecture~\ref{c:CHH} in the affirmative 
with the following theorem.
% -------------------------------------------------------------
\begin{theorem}\label{t:main}
Let $d$ be a positive integer and set $k=2^{2d+3}$.
If $\GVE$ is a graph of
diameter at most~$d$ and connectivity at least~$k$,
then the pebbling number of $\bfG$ is $|V|$.
\end{theorem}
% -------------------------------------------------------------
\begin{proof}
Let $d$ be a positive integer and set $k=2^{2d+3}$. 
Let $\GVE$ be any graph with diameter at most $d$ and 
connectivity at least $k$.  Then let $\phi$ be any pebbling of 
size $|V|$ on $\bfG$.  We assume that there is a
vertex $z_0\in V$ so that $\phi$ does not have the
potential to pebble $z_0$ and argue to a contradiction.

We begin by defining a partition of the vertex set of $\bfG$ by
setting $V=Z\cup U\cup B$, where 
\begin{enumerate}
\item $Z=\{z\in V:\phi(z)=0\}$,
\item $U=\{u\in V:\phi(u)=1\}$, and
\item $B=\{b\in V:\phi(b)>1\}$.
\end{enumerate}
In the remainder of the argument, we use the natural convention
that vertices of $Z$ (\textit{zeroes}) will be denoted by
the letter $z$ (perhaps with subscripts or primes appended).
Similarly, elements of $U$ (\textit{units})
will be denoted by the letter $u$, while elements of $B$ 
(\textit{bigs}) will be denoted by the letter $b$.  Whenever
we want to make a statement about an arbitrary vertex of the
graph, we will use the letter $v$.

Next, we observe that since
the size of $\phi$ is $|V|$, we must have
\[
\sum_{b\in B}\phi(b)=|B|+|Z|.
\]
Of course the vertex $z_0$ belongs to $Z$, so $Z\neq\emptyset$.
Thus $B\neq\emptyset$.  Note that there are no edges from
$z_0$ to vertices in $B$.   Let $m=|B|$.  Then let $\omega=
\sum_{b\in B}\phi(b)/m$ denote the average number of pebbles
assigned by $\phi$ to vertices in $B$.  
It follows that $|Z|=m(\omega - 1)$.  

Before proceeding with the proof, we pause to make a
few elementary observations about properties which
the pebbling $\phi$ must satisfy.  As noted previously,
the vertex $z_0$ has no neighbors in $B$.  In fact,
more can be said.  There cannot be a path $P$ beginning
at $z_0$ and ending at a point in $B$ with all interior
points of $P$ belonging to $B\cup U$.  This follows
from the fact that the existence of such a path would allow
us to shift a pebble from the endpoint of $P$ which
belongs to $B$ along the path until it rests on $z_0$.  

Also, as noted previously, we know that $\phi(v)<2^d$ for every 
$v\in V$.  But again we can say more.  

\medskip
\noindent\textbf{Claim 1.}
Let $v\in V$ and let
$\mathcal{E}$ be a family of paths such that
\begin{enumerate}
\item Each path in $\mathcal{E}$ begins at $v$ and ends
at a point of $B$.
\item For each $b\in B$, at most $\lfloor\phi(b)/2\rfloor$
paths in $\mathcal{E}$ end at $b$.
\item No two paths in $\mathcal{E}$ have any
points in common, apart from $v$ and, if they end at the same point,
their common endpoint.
\item All interior points of paths in $\mathcal{E}$ belong to $B\cup U$.
\end{enumerate}
Then $|\mathcal{E}|<2^d$.

\medskip\noindent\textit{Proof.}
Let $b\in B$.  Then let $\mathcal{E}_b$ denote the
set of all paths in $\mathcal{E}$ which end at
$b$.  We know that $|\mathcal{E}_b|\le\lfloor\phi(b)/2\rfloor$.
It follows that we may shift $|\mathcal{E}_b|$ pebbles from
$b$ to $v$, one along each path in $\mathcal{E}_b$.
Since the paths in $\mathcal{E}$ have no interior points
in common, it follows that a descendant of $\phi$
places  $|\mathcal{E}|$ pebbles on $v$.  This requires
$|\mathcal{E}|<2^d$, as claimed.\hfill$\triangle$

\medskip
Now let $v\in V$ and consider the subgraph $\bfH_v$
of $\bfG$ induced by $\{v\}\cup U\cup B$.  We modify
$\bfH_v$ into a new graph $\bfH'_v$ as follows.   For 
each vertex $b\in B-\{v\}$, we replace $b$ by
an independent set $A_b$ of cardinality $\lfloor\phi(b)/2\rfloor$.
Furthermore, if $b$ is adjacent to a vertex $v'$ in $\bfH_v$,
then every vertex of $A_b$ is adjacent to $v'$ in $\bfH'_v$,
and if $b_1,b_2\in B-\{v\}$ and $b_1$ is adjacent to
$b_2$ in $\bfH_v$, then every vertex of $A_{b_1}$ is adjacent
to every vertex of $A_{b_2}$ in $\bfH'_v$.  We let $B'_v=
\cup\bigl\{A_b:b\in B-\{ v\}\bigr\}$.  Then add a new
vertex $\hat{v}$ with $\hat{v}$ adjacent to all vertices of
$B'_v$ but to no other vertices in $\bfH'_v$.  In particular,
$\hat{v}$ is not adjacent to $v$ in $\bfH'_v$.

We now apply Menger's theorem to the non-adjacent pair
$\{v,\hat{v}\}$, i.e., the minimum number of vertices in
$\bfH'_v$ required to separate $v$ and $\hat{v}$ is equal
to the maximum number of pairwise disjoint paths from $v$ to $\hat{v}$.
Choose a minimum subset $S'_v$ of vertices separating
$v$ from $\hat{v}$ in $\bfH'_v$.  Then
$\{v,\hat{v}\}\cap S'_v=\emptyset$ and
every path in $\bfH'_v$ beginning at 
$v$ and ending at $\hat{v}$ passes through one or 
more points of $S'_v$.  Note that $S'_v\subseteq B'_v\cup U$ for
every $v\in V$.  

\medskip
\noindent\textbf{Claim 2.}
For every $v\in V$, the following statements hold.
\begin{enumerate}
\item Any path in $\bfH'_v$ beginning at $v$ and ending at
a point of $B'_v$ passes through a point of $S'_v$.
\item $|S'_v|<2^d$.
\item If $b\in B$, $b\neq v$ and $A_b\cap S'_v\neq\emptyset$,
then $A_b\subseteq S'_v$.
\end{enumerate}  

\medskip
\noindent
\textit{Proof.}  Let $v\in  V$.  The fact that statement~1 holds
is an immediate consequence of Menger's theorem and the definition
of the graph $\bfH'_v$.  We now show that statement~2 holds.

Let $t=|S'_v|$. Then there is a family $\mathcal{G}$ of 
$t$ paths in $\bfH'_v$ each beginning at $v$ and ending at
$\hat{v}$ with any two of these paths having no point in common other 
than $v$ and $\hat{v}$.   For each path $P$ in $\mathcal{G}$,
start at $v$ and travel along $P$ towards $\hat{v}$.
Then let $b(P)$ be the first point on $P$ which belongs to $B'_v$.
Clearly, there must be such a point since the next to last point of
$P$ belongs to $B'_v$.  Note that any point of $P$ between $v$ and
$b(P)$ on $P$ belongs to $U$.  
Then let $P'$ denote the initial segment of 
$P$ beginning at $v$ and ending at $b(P)$,  
and let $\mathcal{G}'=\{P':P\in\mathcal{G}\}$. 
Note that by Menger's theorem,
each path in $\mathcal{G}'$ contains a unique point
of $S'_v$.

For each $b\in B-\{v\}$, let $\mathcal{G}'_b$
denote the set of all paths in $\mathcal{G}'$ which end at
a point of $A_b$, and let $B_v$  consist of those points
of $B-\{v\}$ for which $\mathcal{G}'_b\neq\emptyset$.
Clearly, $|\mathcal{G}'_b|\le|A_b|=
\lfloor\phi(b)/2\rfloor$ for each $b\in B_v$.   

For each $b\in B_v$, let $\mathcal{G}_b$ denote
the family of paths in $\bfH_v$ obtained by replacing the ending
point of each path in $\mathcal{G}'_b$ by $b$.  Evidently, all 
paths in $\mathcal{G}_b$ start
at $v$ and end at $b$.  However, other than starting and
ending points, the paths in $\mathcal{G}_b$ are pairwise disjoint.
Also, if $b_1$ and $b_2$ are distinct
elements of $B_v$, then 
$v$ is the unique common point of a path from $\mathcal{G}_{b_1}$
and a path from $\mathcal{G}_{b_2}$.
From these remarks, it follows from Claim~1 that
$|S'_v|=|\mathcal{G}|=|\mathcal{G}'|<2^d$, as claimed.

Finally, we show that statement~3 holds.  
Let $b\in B-\{v\}$ and suppose that $A_b\cap
S'_v\neq\emptyset$.  Choose a path $P$ in $\mathcal{G}$
having a point $b_1\in A_b\cap S'_v$ in its interior.  Then
$b_1$ is the unique point from $S'_v$ belonging to $P$.  If there
is a point $b_2\in A_b-S'_v$, then the path $\hat{P}$ obtained from
$P$ by appending $b_2$ and $\hat{v}$ onto the initial segment of $P$
beginning at $v$ and ending at the vertex immediately preceding
$b_1$ on $P$ is a path from $v$ to $\hat{v}$ in $\bfH'_v$
which passes through no point of $S'_v$.  This would contradict the 
assumption that $S'_v$ separates $v$ from $\hat{v}$.  
So $A_b\subseteq S'_v$ as claimed.
\hfill$\triangle$

\medskip
For each $v\in V$, let $S_v=(S'_v\cap U)\cup\{b\in B-\{v\}:A_b\subseteq
S'_v\}$.  Note that $v\not\in S_v$.  
Also note that $S_v$ separates 
$v$ from $B-\{v\}$ in $\bfH_v$, i.e., any path in $\bfH_v$ starting at
$v$ and ending at a point in $B-\{v\}$ passes through one or more
points of $S_v$.  From Claim~2, it is clear that 
$|S_v|<2^d$.  However, using the second part of Claim~2, even
more can be said.  

\medskip\noindent\textbf{Claim 3.}
For every $v\in V$,
\[
\phi(v)+\sum_{b\in S_v}\phi(b)<2^{d+2}.
\]

\medskip\noindent\textit{Proof.} Let $v\in V$. Obviously, $\phi(v)<2^d$.
From Claim~2, we know that
\[\sum_{b\in S_v}\lfloor\phi(b)/2\rfloor
= \sum_{b\in S_v}|A_b|\le|S'_v|< 2^d\]
and
\[\sum_{b\in S_v}\phi(b)\le 
2\sum_{b\in S_v}\Bigl(\lfloor\phi(b)/2\rfloor+|S'_v|\Bigr)\le 3(2^d),\]
so that
\[
\phi(v)+\sum_{b\in S_v}\phi(b)<2^d+3(2^d)=2^{d+2}.
\]
\hfill$\triangle$

\medskip
\noindent\textbf{Claim 4.}  There exists a positive integer 
$q$ with $q> m\omega/2^{d+2}$, a $q$-element subset 
$B_0\subset B$, and a labelling $\{b_1,b_2,\dots,b_q\}$ of the
elements of $B_0$ so that for all $i$ and $j$ with $1\le i<j\le q$,
$b_j\not\in S_{b_i}$.  

\medskip
\noindent\textit{Proof.}
We form the subset $B_0$ inductively.  Choose an
arbitrary element of $B$ as $b_1$.  Then
remove from consideration all remaining elements $S_{b_1}\cap B$.  
By Claim~3, the total number of pebbles
assigned by $\phi$ to the elements selected or removed from consideration 
is less than $2^{d+2}$.
Since the total number of pebbles assigned to elements of
$B$ is $m\omega$, we may repeat this procedure $m\omega/2^{d+2}$
times to obtain the desired subset $B_0$.
\hfill$\triangle$

\medskip
The reader should note that if $1\le i < j\le q$, then
we know that $b_j\not\in S_{b_i}$, but
we do not know whether $b_i$ belongs to $S_{b_j}$.
Now let $W= \cup_{b\in B_0}S_b$.  Note that $|W|<
q2^d$.  Also note that $W\cap Z=\emptyset$.

Since the connectivity of $\bfG$ is at least $k$, we know
that for each $b\in B_0$, there are
$k$ paths $P_1(b,z_0), P_2(b,z_0),\ldots, P_k(b,z_0)$, each beginning
at $b$ and ending at $z_0$, with two paths in this family having
no points in common other than $b$ and $z_0$.  
Since $\{b,z_0\}$ is not an edge in $\bfG$,
each path $P_i(b,z_0)$ contains at least one interior point,
and in fact, each $P_i(b,z_0)$ contains at least one interior point 
which belongs to $Z$. 

Since $|S_b|<2^d$ and $k=2^{2d+3}$, we may assume that
the paths have been labelled so that
for each $i=1,2,\dots,k-2^d$, the path $P_i(b,z_0)$ does not
contain a point of $S_b$. 

Now let $b\in B_0$ and let $i$ be an integer with $1\le i\le k-2^d$.
Follow the path $P_i(b,z_0)$ beginning at $b$ and let $v_i(b)$ be the first
point (distinct from $b$) on the path which belongs to $W\cup Z$.
We let $Q_i(b)$ be the initial segment
of $P_i(b,z_0)$ beginning at $b$ and ending at $v_i(b)$.
We call $b$ the \textit{root} of this path and
$v_i(b)$ the \textit{terminal point},
and we let $\mathcal{F}=\{Q_i(b):b\in B_0, 1\le i\le k-2^d\}$.
Note that $|\mathcal{F}|=q(k-2^d)$.

Of course, for each $v\in W\cup Z$ and each $b\in B$, there
is at most one integer $i$ with $1\le i\le k-2^d$ for
which the path $Q_i(b)$ has $b$ as its root and
$v$ as its terminal point.  However, if $b,b'\in B$, there may
exist integers $i,i'\in\{1,2,\dots,k-2^d\}$ so that
$Q_i(b)$ and $Q_{i'}(b')$ both have $v$ as their terminal point;
but when $b$ and $b'$ both belong to $B_0$, we can say much more.

\medskip
\noindent\textbf{Claim 5.} 
Let $v\in W\cup Z$ and let $b$ and $b'$ be distinct 
points of $B_0$.  If $i,i'\in\{1,2,\dots,k-2^d\}$ and
$Q_i(b)$ and $Q_{i'}(b')$ both have $v$ as their terminal point,
then $Q_i(b)$ and $Q_{i'}(b')$ have no point in common other
than $v$. 

\medskip
\noindent
\textit{Proof.}  
Suppose that $F =Q_i(b)$ and $F'=Q_{i'}(b')$ have 
a common point distinct from $v$.  Then there exists a path $F''$ 
from $b$ to $b'$ with $F''\subseteq (F\cup F')-\{v\}$.  
Furthermore, $F''\subseteq B\cup U$.
Thus $F''\cap S_{b}\neq\emptyset\neq F''\cap S_{b'}$.  However,
$F''\cap F$ does not contain a point from $S_{b}$, and
$F''\cap F'$ does not contain a point from $S_{b'}$.
Therefore $F''\cap F\cap S_{b'}\neq\emptyset$ and 
$F''\cap F'\cap S_{b}\neq\emptyset$.  

Without loss of generality, we may assume that $b$ was chosen
before $b'$ in the construction of $B_0$ described in
Claim~3.  Then $b'\not\in S_b$.  It follows that $F''\cap F'\cap S_b$
contains a point of $S_b$ which is distinct from $b'$.
However, this contradicts the hypothesis that $v$ is the
terminal point of $F'$.
\hfill$\triangle$

\medskip
We pause to point out that the argument in the preceding claim
works only in one direction, as it may happen that
$b\in S_{b'}$.

\medskip
We are now ready to complete the proof.  We observe that
$|W\cup Z| <q2^d+m(\omega-1)<q(2^{d+2}+2^d)$ and 
$|\mathcal{F}|=q(k-2^d)$.  Since $|\mathcal{F}|/|W\cup Z|>2^d-1$,
it follows that there is a vertex $v\in W\cup Z$ and a subfamily
$\mathcal{E}\subset\mathcal{F}$, with $|\mathcal{E}|=2^d$, so
that every path in $\mathcal{E}$ has $v$ as its terminal point.
However, from Claim~5,
any two paths from $\mathcal{E}$ have no 
point in common other than $v$, so the existence
of $\mathcal{E}$ is now seen to be a contradiction to
Claim~1.
\end{proof}
% =============================================================
\section{Threshold}\label{s:thresh}
% =============================================================
The notion that graphs with very few edges tend to have large pebbling
number and graphs with very many edges tend to have small pebbling
number can be made precise as follows.
Let $\gnp$ be the random graph model in which each of the $n\choose 2$ 
possible edges of a random graph having $n$ vertices appears independently 
with probability $p$.
For functions $f$ and $g$ on the natural numbers we write that
$f\ll g$ (or $g\gg f$) when $f/g\rightarrow 0$ as $n\rightarrow\infty$.
Let $o(g)=\{f\ |\ f\ll g\}$ and define $O(g)$ (resp., $\Omega(g)$) to be 
the set of functions $f$ for which there are constants $c, N$ such that
$f(n)\le cg(n)$ (resp., $f(n)\ge cg(n)$) whenever $n>N$.
Finally, let $\Theta(g)=O(g)\cap\Omega(g)$.

Let $\cgQ$ be a property of graphs and consider the probability
$Pr(\cgQ)$ that the random graph $\gnp$ has $\cgQ$.
For large $p$ it may be that $Pr(\cgQ)\rightarrow 1$ as 
$n\rightarrow\infty$, and for small $p$ it may be that 
$Pr(\cgQ)\rightarrow 0$ as $n\rightarrow\infty$.
More precisely, define the \textit{threshold of} $\cgQ$, $th(\cgQ)$, 
to be the set of functions $t$ for which $p\gg t$ implies that
$Pr(\cgQ)\rightarrow 1$ as $n\rightarrow\infty$, and $p\ll t$ implies that
$Pr(\cgQ)\rightarrow 0$ as $n\rightarrow\infty$, if this set is nonempty.

It is not clear that such thresholds exist for arbitrary $\cgQ$.
However, we observe that Class 0 is a monotone property (adding edges
to a Class 0 graph maintains the property), and a theorem of Bollob\'as
and Thomason \cite{BT} states that $th(\cgQ)$ exists for every 
monotone $\cgQ$. 
It is well known \cite{ER1} that $th({\rm connected})=\Theta(\lgn /n)$, 
and since connectedness is required for Class 0, we see that 
$th({\rm Class\ 0})\subseteq\Omega(\lgn /n)$.
In \cite{CHH} it is noted that $\gn$ is Class 0 with probability tending to 1.
It is straightforward to extend their argument to show that,
for fixed $p$, $\gnp$ is Class 0 with probability tending to 1.
Here we prove the following theorem.
% -------------------------------------------------------------
\begin{theorem}\label{t:cor}
For all $d>0$, $th({\rm Class\ }0)\subseteq o((n\lgn)^{1/d}/n)$.
\end{theorem}
% -------------------------------------------------------------
\begin{proof}
We prove the equivalent statement that
$th({\rm Class\ }0)\subseteq O((n\lgn)^{1/d}/n)$ for all $d>0$.
It is proven in \cite{B} that
$th({\rm diameter} \le d)\subseteq\Omega((n\lgn)^{1/d}/n)$,
and in \cite{ER2} that 
$th({\rm connectivity} \ge k)\subseteq\Omega((\lgn +k\lg{\lgn})/n)$.
Hence, for any fixed $d$ and $k$ with $k\ge 2^{2d+3}$, and for any
$p\gg (n\lgn)^{1/d}/n$, the probability that $\gnp$ has diameter at most
$d$ and connectivity at least $k$ tends to 1.
Therefore the probability that $\gnp$ is Class 0 tends to 1.
\end{proof}

% =============================================================
\section*{Acknowledgements}
% =============================================================
The authors are grateful for the many useful comments provided by the
referees.
% -------------------------------------------------------------


\begin{thebibliography}{99}
% -------------------------------------------------------------
\bibitem{BT}
B. Bollob\'as and A. G. Thomason, 
Threshold functions,
\textit{Combinatorica} \textbf{7} (1987), 35--38.
\bibitem{B}
Y. D. Burtin,
On extreme metric characteristics of a random graph II:
Limit distributions,
\textit{Theory Probab. Appl.} \textbf{20}, 83--101.
\bibitem{CHH}
T. A. Clarke, R. A. Hochberg and G. H. Hurlbert,
Pebbling in diameter two graphs and products of paths,
\textit{J. Graph Theory} \textbf{25} (1997), 119--128.
\bibitem{ER1}
P. Erd\H os and A. R\'enyi,
On random graphs I,
\textit{Publ. Math. Debrecen} \textbf{6} (1959), 290--297.
\bibitem{ER2}
P. Erd\H os and A. R\'enyi,
On the strength of connectedness of a random graph,
\textit{Acta Math. Sci. Hung.} \textbf{12} (1961), 261--267.
% -------------------------------------------------------------
\end{thebibliography}
\end{document}